\documentclass[a4paper,reqno, 11pt]{amsart}
\author{Julia Brandes}
\thanks{Partially supported by an EPSRC DTA studentship}
\title{A note on $p$-adic solubility for forms in many variables}

\address{Mathematisches Institut, Bunsenstr. 3--5, 37073 G{\"o}ttingen, Germany}
\email{jbrande@uni-math.gwdg.de}

\usepackage[T1]{fontenc}
\usepackage[latin1]{inputenc}
\usepackage{latexsym}
\usepackage[arrow, matrix, curve]{xy}
\usepackage[dvips]{graphicx}
\usepackage{epsfig}
\usepackage{bm}
\usepackage{amssymb}
\usepackage{amsmath}
\usepackage{verbatim}
\usepackage{amsthm}
\usepackage{enumerate}
\usepackage{mathrsfs}
\usepackage[hmargin=3.5cm,vmargin=3.5cm]{geometry}

\newtheorem{thm}{Theorem}
\newtheorem{lem}{Lemma}

\newtheorem{cor}{Corollary}

\numberwithin{equation}{section}
\numberwithin{thm}{section}
\numberwithin{lem}{section}	
\numberwithin{cor}{section}
\theoremstyle{definition}

\def\B#1{\mathbf{#1}}

\def\F#1{\mathfrak{#1}}
\def\dsum#1#2{\sum_{\substack{{#1}\\{#2}}}}
\def\ol#1{\overline{\B{#1}}}

\def\U#1{\underline{\B{#1}}}

\DeclareMathOperator{\rk}{rank}
\DeclareMathOperator{\card}{Card}
\DeclareMathOperator{\Id}{Id}
\DeclareMathOperator{\ord}{ord}

\usepackage{hyperref}

\begin{document}

\begin{abstract}

    By adopting a new approach to the analysis of the density of $p$-adic solutions arising in applications of the circle method, we show that under modest conditions the existence of non-trivial $p$-adic solutions suffices to establish positivity of the singular series. This improves on earlier approaches due to Davenport, Schmidt and others, and allows us to establish an asymptotic formula for the number of simultaneous zeros of non-singular pairs of cubic forms in at least $131$ variables. As a by-product, we obtain a version of Hensel's Lemma for linear spaces. 
\end{abstract}

\maketitle

\section{Introduction}

Thanks to the groundbreaking work of Davenport \cite{dav32} and Birch \cite{birch} and later Schmidt \cite{schmidt85}, the treatment of forms in many variables is now one of the most classical applications of the circle method. However, as the generic outcome of the method is only a local-global principle, additional information on the density of solutions in the local fields $\mathbb R$ and $\mathbb Q_p$ is required in order to derive a proper asymptotic formula with a non-vanishing main term for the number of solutions over $\mathbb Q$. Unfortunately, the available methods for finding non-trivial $p$-adic solutions do not address their quantitative behaviour, and whilst it is not impossible to obtain information on the solution density from our knowledge about non-trivial $p$-adic solubility, the technical devices which are commonly applied in order to do so (e.g. \cite{cubIV}) are complicated and wasteful. 
In this note we address this problem and present a more efficient geometric approach for making the transition from non-trivial to non-singular $p$-adic points at the cost of inflating the number of variables by no more than what is required by the geometry of the problem. 

Let $F^{(1)}, \ldots, F^{(R)} \in \mathbb Z[x_1, \ldots, x_s]$ be homogeneous polynomials of common degree $d$. We follow Dietmann~\cite{rd-weyl} und Schindler \cite{damaris-weyl} in defining 
$$\F V_{\B a}(\B F) = \left\{\B x \in \mathbb C^s: \nabla\big(a_1F^{(1)}(\B x) + \dots + a_RF^{(R)}(\B x)\big) = 0 \right\}$$ 
for $\B a \in \mathbb Z^R \setminus \{\bm 0\}$. Furthermore, let $\F S(\B F)$ denote the singular locus of the variety given by $\B F(\B x)=\bm 0$. If $N(P)$ denotes the number of integral points $\B x$ with $|x_i| \leq P$ for all $1 \leq i \leq s$ satisfying the equations $F^{(1)}(\B x) = \ldots = F^{(R)}(\B x) = 0$, we write $H_d(R)$ for the least integer $H$ with the property that every system $\B F$ as above satisfying $s - \max_{\B a \ne \bm 0}\dim \F V_{\B a}(\B F) \geq H $ obeys an asymptotic formula of the shape 
\begin{align}\label{formula}
 N(P) = (c + o(1)) P^{s-Rd},
\end{align}
where $c$ is a non-negative constant encoding the solution densities over the local fields $\mathbb R$ and $\mathbb Q_p$.  
Similarly, let $\gamma^p_d(R)$ denote the least integer $\gamma$ such that any system of $R$ forms of degree $d$ in $s \geq \gamma$ variables has  a non-trivial simultaneous zero in $\mathbb Q_p$, and write $ \gamma_d^*(R) = \max_{p} \gamma^p_d(R)$. We stress here that $\gamma_d^*(R)$ is defined in terms of $s$ itself and does not depend on singularities of any kind.

\begin{thm}\label{thm:asymp}
    Let $F^{(1)}, \ldots, F^{(R)}\in \mathbb Z[x_1, \ldots, x_s]$ be forms of degree $d \ge 3$, and suppose that 
    \begin{align*}
      s - \max_{\B a  \in \mathbb Z^R \setminus \{\bm 0\}}\dim \F V_{\B a}(\B F)\geq H_d(R) \quad \text{ and } \quad s - \dim\F S(\B F) \ge \gamma_d^*(R).
    \end{align*}
    Then \eqref{formula} is satisfied with a positive constant $c$, provided the system has a non-singular real solution. In particular, this will be the case if $d $ is odd.
\end{thm}
Notice that Birch's Theorem~\cite{birch} in conjunction with the results of \cite{rd-weyl,damaris-weyl} provides the bound
\begin{equation*}
    H_d(R) \leq  2^{d-1} (d-1)R(R+1) +1,
\end{equation*} 
and in the case when $d$ is odd the density of solutions in $\mathbb R$ can be shown to be positive by an argument due to Schmidt (see \cite[\S 2]{cubIV} and \cite[Lemma~2 and \S 11]{schmidtquad}) which also requires $s-\max_{\B a \ne \bm 0}\dim \F V_{\B a}(\B F) \geq H_d(R)$. The significance of Theorem~\ref{thm:asymp} lies therefore in the implied bound $s - \dim \F S(\B F)\geq \gamma_d^*(R)$ to guarantee a positive solution density in all fields $\mathbb Q_p$. This bound clearly supersedes the previous treatment due to Schmidt (see  \cite[\S 3]{cubIV}), which in our setting gives the condition 
\begin{equation*}
    s - \max_{\B a  \in \mathbb Z^R \setminus \{\bm 0\}}\dim \F V_{\B a}(\B F)>  2^{d-1} (d-1)R \gamma_d^*(R).
\end{equation*}
In comparison, our bound is much simpler, such as one would na\"{i}vely expect, and indeed the effect of our new result is that from a philosophical point of view the non-diagonal problem is now on an equal footing with Waring-type situations. 
A similar result (see Theorem \ref{thm:asymp-m} below) is available for the number of $m$-dimensional linear spaces of height $P$ that are contained in the complete intersection defined by $\B F= \bm 0$. 

Theorem~\ref{thm:asymp} can be made explicit by consulting the literature. 
\begin{cor}\label{cor:birch}
    \begin{enumerate}
	\item Let $\B F$ be as in Theorem~\ref{thm:asymp}. Then \eqref{formula} holds for a positive constant $c$, provided that ${s - \dim \F S(\B F)\geq (Rd^2)^{2^{d-1}}}$.
	\item Let $F$ and $G$ be cubic forms in $s$ variables, then the number of solutions to $F(\B x)=G(\B x)=0$ satisfies \eqref{formula} with a positive constant $c$ as soon as $s - \dim \F S(F, G)\geq 131$.
    \end{enumerate}
\end{cor}
The first statement is immediate from Wooley's bound in \cite[Cor. 1.1]{tdw-local} and the fact that $ (Rd^2)^{2^{d-1}} > 2^{d-1} (d-1)R(R+1)+R-1$ for all positive parameters $d$ and $R$, where we used the fact that 
\begin{align*}
    \max_{\B a  \in \mathbb Z^R \setminus \{\bm 0\}}\dim \F V_{\B a}(\B F) \le \dim\{\B x: \rk(\partial_i F^{\rho}(\B x))_{i, \rho} \le R-1\} \le \dim \F S(\B F) + R-1.
\end{align*}
The second bound follows from a recent result by Dumke \cite{dumke} which establishes ${\gamma^*_3(2) \leq 131}$. 
Previously, the best known bound was a result by Dietmann and Wooley \cite{rd-tdw}, who by a different approach showed that any two cubic forms, not necessarily non-singular, have a simultaneous zero if the number of variables is at least $827$. Whilst this is more general than our result in that it does not require non-singularity, in the generic case Corollary \ref{cor:birch} gives a stronger result both in terms of quantitative behaviour and the number of variables required. Note that in the case of a single non-singular cubic form Heath-Brown \cite{hb10} was able to show that as few as ten variables suffice, and this is best possible. \\

I would like to express my gratitude towards my PhD supervisor Trevor Wooley for his keen insight and constant encouragement, and to my examiners Tim Browning and Rainer Dietmann as well as a pair of anonymous referees for valuable comments. This work forms part of the author's PhD thesis.\\

\section{Setup}

We will phrase the entire argument so as to accommodate higher-dimensional linear spaces as well as points. For this purpose it is necessary to understand expressions of the shape
\begin{equation}\label{linspace-gen}
    F^{(\rho)}(t_1 \B x_1+ \ldots + t_m \B x_m).
\end{equation}
Suppose the forms $F^{(1)}, \ldots, F^{(R)} \in \mathbb Z[x_1, \ldots, x_s]$ are given by
\begin{equation*}
    F^{(\rho)}(\B x)=\sum_{\B i \in \{1, \ldots,s\}^d} c^{(\rho)}_{\B i} x_{i_1}  \cdot \ldots \cdot  x_{i_d} \quad (1 \leq \rho \leq R)
\end{equation*}
with symmetric coefficients $c^{(\rho)}_{\B i} \in \mathbb Z / d!$. Then we can define the symmetric multilinear forms $\Phi^{(\rho)}$ associated to $F^{(\rho)}$ as
\begin{equation*}
    \Phi^{(\rho)}\big(\B x^{(1)}, \ldots ,\B x^{(d)}\big) = \sum_{\B i \in \{1, \ldots, s\}^d} c^{(\rho)}_{\B i} x_{i_1}^{(1)}  \cdot \ldots \cdot  x_{i_d}^{(d)},
\end{equation*}
such that $F^{(\rho)}(\B x) = \Phi^{(\rho)}(\B x, \ldots, \B x)$ for every $\rho$.
Write $J$ for the set of multi-indices
\begin{equation*}
    J=\lbrace (j_1, j_2, \ldots, j_d) \in \{1, 2, \ldots, m\}^d : j_1 \leq j_2 \leq \ldots \leq j_d\rbrace;
\end{equation*}
the number of these is 
\begin{equation*}
    \card(J)  = \binom{d-1+m}{d} = r.
\end{equation*}
 By means of the Multinomial Theorem, equation \eqref{linspace-gen} can be written as
\begin{equation}\label{expanded}
    F^{(\rho)}\left(t_1 \B x_1+ \ldots + t_m \B x_m\right)
    =\sum_{ \B j \in J} A(\B j) t_{j_1}t_{j_2}  \cdot \ldots \cdot t_{j_d} \Phi^{(\rho)}(\B x_{j_1}, \B x_{j_2},\ldots, \B x_{j_d})
\end{equation}
with suitable combinatorial factors $A(\B j)$ that take account of the multiplicity of each term.

By writing $ \left(\B x_1, \ldots, \B x_m\right) =\ol x$ and
\begin{equation*}
    \Phi^{(\rho)}_{\B j}(\ol x) = A(\B j) \Phi^{(\rho)}(\B x_{j_1}, \ldots \B x_{j_d})  \in \mathbb Z[\B x_1, \ldots, \B x_m] \qquad (1 \leq \rho \leq R, \, \B j \in J),
\end{equation*}
it follows from equation \eqref{expanded} that every linear space of dimension $\le m$ contained in the complete intersection defined by $\B F = \bm 0$ corresponds to a solution of the system of equations
\begin{equation}\label{eq:sys}
    \Phi^{(\rho)}_{\B j}(\ol x) = 0 ,\qquad \B j \in J, \, 1 \leq \rho \leq R.
\end{equation} 
Finally, if $\Gamma_m(p^l)$ denotes the number of solutions $\B x_1, \ldots, \B x_m \in (\mathbb Z/p^l \mathbb Z)^s$ of the simultaneous congruences
\begin{equation*}
    F^{(\rho)}(\B x_1 t_1 + \ldots + \B x_m t_m) \equiv 0 \pmod{p^l}, \quad 1 \leq \rho \leq R,
\end{equation*}
identically in $t_1, \ldots, t_m$, then by \eqref{expanded} this can be written as
\begin{equation*}
    \Gamma_m(p^l) = \card \left\lbrace \ol x \in (\mathbb Z/p^l\mathbb Z)^{ms} :  \Phi^{(\rho)}_{\B j}(\ol x) \equiv 0 \pmod{p^l}  \quad \hbox{ for all } \B j ,\rho  \right\rbrace. 
\end{equation*}
We are now equipped with all the necessary notation to embark on the argument.

\section{Hensel's Lemma~for linear spaces}\label{sec:hensel}

The proof of Theorem~\ref{thm:asymp} depends crucially on Hensel's Lemma, and as we are aiming to phrase our result in terms of linear spaces, we need a suitable multi-dimensional version of it. As is usual with arguments in the spirit of Hensel's Lemma, this involves finding a non-singular approximate solution which can then be lifted to higher moduli and thus generate solutions in~$\mathbb Q_p$. \\ %

Suppose $\ol a \in (\mathbb Z_p)^{ms}$ is a solution to \eqref{eq:sys}.
An application of Taylor's Theorem~to each of the equations in \eqref{eq:sys} yields
\begin{align}
    \Phi^{(\rho)}_{\B j} (\ol a + p^{h}\ol x) & =\Phi^{(\rho)}_{\B j} (\B a_{j_1} + p^{h}\B x_{ j_1}, \B a_{j_2} + p^{h}\B x_{ j_2}, \ldots, \B a_{j_d} + p^{h}\B x_{ j_d}) \nonumber\\
    &= \Phi^{(\rho)}_{\B j}(\ol a) + p^{h}L^{(\rho)}_{\B j}(\ol a) \cdot \ol x + R^{(\rho)}_{\B j}(\ol a) (p^{h} \ol x), \label{eq:taylor}
\end{align}
where $\Phi^{(\rho)}_{\B j}(\ol a) = 0$ by assumption, $R^{(\rho)}_{\B j}(\ol a) (\ol x)$ is a polynomial in $\ol x$ without constant or linear terms, and $L^{(\rho)}_{\B j}(\ol a)$ is defined as follows. If $\B e_n$ is the $n$-th unit vector, we write
\begin{equation*}
    B^{(\rho)}_n(\B x^{(1)}, \ldots ,\B x^{(d-1)}) = \Phi^{(\rho)}( \B e_n, \B x^{(1)}, \ldots ,\B x^{(d-1)}), \quad 1 \leq \rho \leq R,
\end{equation*}
so that
\begin{equation*}
    \Phi^{(\rho)}(\B x, \B a_1, \ldots, \B a_{d-1}) = \sum_{n=1}^s B_n^{(\rho)}(\B a_1, \ldots, \B a_{d-1}) x_n.
\end{equation*}
Then $L^{(\rho)}_{\B j}(\ol a)$ is defined via the relation
\begin{align*}
    L^{(\rho)}_{\B j}(\ol a) \cdot \ol x &= A(\B j) \sum_{k=1}^d  \Phi^{(\rho)}(\B a_{j_1}, \ldots, \B a_{j_{k-1}}, \B x_{j_k}, \B a_{j_{k+1}}, \ldots, \B a_{j_d}) \\
    &= A(\B j)\sum_{k=1}^d \sum_{n=1}^s B^{(\rho)}_n(\B a_{j_1}, \ldots, \B a_{j_{k-1}}, \B a_{j_{k+1}}, \ldots, \B a_{j_d}) x_{j_{k}, n}   \\
    &= A(\B j)\sum_{i=1}^m \sum_{n=1}^s \left[  \sum_{k=1}^d B^{(\rho)}_n(\B a_{j_1}, \ldots, \B a_{j_{k-1}}, \B a_{j_{k+1}}, \ldots, \B a_{j_d}) \delta_{j_k, i} \right] x_{i,n},
\end{align*}
where $\delta_{i, j} $ denotes the Kronecker delta. Writing $L^{(\rho)}_{\B j, n, i}(\ol a)$ for the term in the square brackets, one has the $(Rr\times ms)$-matrix
\begin{equation*}
    \mathscr{L}(\ol a) = \left[ L^{(\rho)}_{\B j, n, i}(\ol a) \right]_{\B j, \rho; n, i} = \left(\frac{\partial \Phi_{\B j}^{\rho}(\ol a)}{\partial \B a_i} \right)_{\B j, \rho; i}.
\end{equation*}
A solution $\ol a$ of the system \eqref{eq:sys} is singular if it satisfies $\rk (\mathscr{L}(\ol a)) \leq Rr-1$, and we may define the $m$-dimensional singular locus of the system $\B F$ as the affine variety
\begin{equation*}
    \F S_m(\B F) = \left\{\ol x: \, \Phi_{\B j}^{(\rho)}(\ol x) = 0 \quad (\B j \in J, 1 \le \rho \le R) \quad \text{ and} \quad \rk (\mathscr{L}(\ol x)) \leq Rr-1 \right\}.
\end{equation*}
Observe that in particular a point $\ol x = (\B x_1, \dots, \B x_m)$ is singular whenever the vectors $\B x_1, \dots, \B x_m$ fail to be linearly independent. 

For an integral matrix $A \in \mathbb Z^{n \times m}$ with $n \leq m$ let $\ord_p(A)$ denote the least integer $h$ such that at least one of the $(n\times n)$-minors of $A$ has determinant not divisible by $p^h$, if such a number exists, and $\ord_p(A)=\infty$ otherwise. If the matrix~$A$ depends on an integral vector $\B x$ we will write $\ord_p(A) = \displaystyle{\min_{\B x } } \ord_p(A(\B x))$. This means that when the system \eqref{eq:sys} is non-singular in $\mathbb Q_p$, it is in some sense non-singular over $\mathbb Z/p^h \mathbb Z$ for all powers $h$ satisfying $h \geq \ord_p(\mathscr{L})$.\\

Suppose $\ord_p(\mathscr L) = \sigma < \infty$, so the system is non-singular, and assume that $ms \ge Rr$. Then for some $\ol a \in \mathbb Z^{ms}$ the matrix $\mathscr L(\ol a)$ has an $(Rr \times Rr)$-minor $\mathscr L_0(\ol a)$ with the property
\begin{align*}
    \left| \det \mathscr L_0(\ol a) \right|_p &= p^{1-\sigma}
\end{align*}
where $| \cdot |_p$ denotes the usual $p$-adic absolute value; we can assume without loss of generality that $\mathscr L_0(\ol a)$ is the first minor.
Write $M(\sigma,\nu)$ for the number of $\ol a \in (\mathbb Z/ p^{2\sigma-1+\nu}\mathbb Z)^{ms}$ that satisfy
\begin{align}
    \Phi^{(\rho)}_{\B j}(\ol a) &\equiv 0 \pmod{p^{2\sigma-1+\nu}} \qquad \mbox{ for all } \B j \in J, 1 \leq \rho \leq R,  \label{eq:eq-p}\\
    \left| \det \mathscr L_0(\ol a) \right|_p &= p^{1-\sigma}\label{eq:det-p},
\end{align}
and additionally have the property that the vectors $(p^{\sigma-1} \ol a_1, \ol a_2)$ are mutually distinct modulo $p^{2\sigma-1+\nu}$, where we wrote $\ol a_1$ to denote the first $Rr$ entries of $\ol a$ and $\ol a_2$ for the remaining $ms-Rr$ values.
Then $M(\sigma,\nu)$ can be bounded below.

\begin{lem}\label{lem:lift}
    Suppose $ms\geq Rr$. For any $\nu \geq 0$ we have
    \begin{equation*}
        M(\sigma,\nu) \geq p^{(ms-Rr)\nu}M(\sigma,0).
    \end{equation*}
\end{lem}

\begin{proof}
    The proof is an amalgam of the arguments of \cite[Lemma~17.1]{dav} and \cite[Prop.~5.20]{greenberg}. For $M(\sigma,0)=0$ the result is immediate, so it suffices to consider the case ${M(\sigma,0)>0}$. Also, the lemma is trivially true for $\nu = 0$. We can therefore proceed by induction and investigate $ M(\sigma,\nu+1)$ under the assumption that for some given $\nu \geq 0$ the statement is true. Let $\ol a$ be one of the solutions to \eqref{eq:eq-p} and \eqref{eq:det-p} counted by $M(\sigma,\nu)$. As in \eqref{eq:taylor}, Taylor's Theorem~yields
    \begin{equation}\label{eq:cong}
        \Phi^{(\rho)}_{\B j}(\ol a + p^{\sigma+\nu}\ol x) \equiv \Phi^{(\rho)}_{\B j}(\ol a) + p^{\sigma+\nu} L^{(\rho)}_{\B j}(\ol a) \cdot \ol x \pmod{p^{2\sigma+\nu}}, \quad  \B j \in J, 1 \leq \rho \leq R,
    \end{equation}
    and by \eqref{eq:eq-p} we can  write $\Phi^{(\rho)}_{\B j}(\ol a) = p^{2\sigma-1+\nu} \varphi_{\B j}^{(\rho)}$ for some $\varphi^{(\rho)}_{\B j} \in \mathbb Z$. Also, by \eqref{eq:det-p} there is a unimodular matrix $\mathscr N \in \mathbb Z^{Rr \times Rr}$ such that
    \begin{equation*}
        \mathscr N \mathscr L_0(\ol a) = \det (\mathscr L_0(\ol a)) \Id_{Rr} = \beta p^{\sigma-1}\Id_{Rr}
    \end{equation*}
    for some $\beta$ coprime to $p$, so the system of congruences
    \begin{equation*}
         \Phi^{(\rho)}_{\B j}(\ol a) + p^{\sigma+\nu} L_{\B j}^{(\rho)}(\ol a) \cdot \ol x \equiv 0\pmod{p^{2\sigma+\nu}}, \quad  \B j \in J, 1 \leq \rho \leq R
    \end{equation*}
    is equivalent to
    \begin{equation}\label{mat1}
        p^{2\sigma-1+\nu} \mathscr N \bm{\varphi}  + p^{\sigma+\nu}\mathscr N \mathscr L(\ol a) \cdot \ol x \equiv \bm 0 \pmod{p^{2\sigma+\nu}}.
    \end{equation}
    The matrix $\mathscr N \mathscr L(\ol a)$ can be written in block form
    \begin{equation*}
        \mathscr N \mathscr L(\ol a) = \mathscr N \big(\mathscr L_0(\ol a) \big| \mathscr L_1(\ol a) \big) = \big(\beta p^{\sigma-1}\Id_{Rr} \big| \mathscr N \mathscr L_1(\ol a)  \big)
    \end{equation*}
    for some matrix $\mathscr L_1(\ol a)$, so if we restrict ourselves to counting solutions $\ol x$ of the shape
    \begin{equation*}
        x_i = \begin{cases} y_i & 1 \leq i \leq Rr \\ p^{\sigma-1}y_i & Rr+1 \leq i \leq ms ,\end{cases}
    \end{equation*}
    we can write
    \begin{equation*}
        \mathscr N \mathscr L(\ol a) \cdot \ol x = p^{\sigma-1} \mathscr L' \cdot \ol y,
    \end{equation*}
    where $\mathscr L' = \big(\beta \Id_{Rr} \big| \mathscr N \mathscr L_1(\ol a)  \big)$. Thus the number of solutions to \eqref{mat1} is at least as large as the cardinality of the solution set of
    \begin{equation*}
        \mathscr N \bm{\varphi} + \mathscr L' \ol y \equiv \bm 0 \pmod p.
    \end{equation*}
    Since $\mathscr N$, $\bm{\varphi}$ and $\mathscr L'$ are fixed, this can be considered as a non-singular system of $Rr$ linear equations in $ms$ variables and therefore has $p^{ms-Rr}$ solutions. 
    Moreover, one checks that the corresponding solutions $\ol a' = \ol a + p^{\sigma + \nu }\ol x$ of \eqref{eq:cong} have the property that the vectors $(p^{\sigma-1}\ol{a}_1', \ol{a}_2')$ are congruent to $(p^{\sigma-1}\ol a_1, \ol a_2)$ modulo $p^{2\sigma-1+\nu}$ but mutually incongruent modulo $p^{2\sigma+\nu}$. Thus altogether we have
    \begin{equation*}
        M(\sigma,\nu+1) \geq p^{(ms-Rr)} M(\sigma,\nu) \geq p^{(ms-Rr)(\nu+1)} M(\sigma,0)
    \end{equation*}
    as claimed.
\end{proof}

Write
\begin{align*}
		\F X_m(\B F)  &=  \left\lbrace  \ol x \in \mathbb Q^{ms} : \;    \Phi_{\B j}^{(\rho)}(\ol x) = 0 \quad \hbox{ for all } \B j \in J, 1 \leq \rho \leq R \right\rbrace
\end{align*}
for the space of solutions to the system \eqref{eq:sys}. It is clear from our previous considerations that if
$\F X_m(\B F) \setminus \F S_m(\B F) \neq \emptyset$, there exists an integer $\sigma$ such that $M(\sigma,0) > 0$.
We can therefore find a non-singular solution and lift it to a $p$-adic solution by means of Lemma~\ref{lem:lift}. This establishes an alternative version of Hensel's Lemma. Notice also that for $m=1$ one has $r=1$ and hence Lemma~\ref{lem:lift} reproduces standard results such as \cite[Lemma~17.1]{dav}.

\begin{lem}[(Hensel's Lemma for linear spaces)]\label{lem:Hensel}
    Let $F^{(1)}, \ldots, F^{(R)} \in \mathbb{Z}[x_1, \ldots, x_s]$ be forms of degree $d$, and $m$ a positive integer. Furthermore, let $ms \geq rR$ and suppose that one has $\F X_m(\B F) \setminus \F S_m(\B F) \neq \emptyset $. Then there exists a constant $\kappa_p$ such that for any sufficiently large $\nu \in \mathbb N$ one has
    \begin{equation*}
        \Gamma_m(p^{\nu}) \geq \kappa_p p^{\nu(ms-Rr)}.
    \end{equation*}
    In particular, the variety defined by $\B F= \bm 0$ contains a non-singular $p$-adic linear space of dimension $m$.
\end{lem}
Note that Lemma~\ref{lem:lift} allows us to take $\kappa_p = p^{(1-2 \sigma)(ms-Rr)}$ for all $\nu \geq 2\sigma-1$, where $\sigma = \ord_p(\mathscr L)$ is a finite parameter depending only on $p$, the dimension $m$ and the system $F^{(1)}, \ldots, F^{(\rho)}$. Observe also that the resulting linear space will be exactly of dimension $m$, as this is a consequence of non-singularity.

\section{Understanding the singularities}

In its present shape, our new version of Hensel's Lemma~is only of limited use, as it is not well adapted to our general state of knowledge. In fact, while there are many results available about the least number of variables required to guarantee the existence of non-trivial $p$-adic solutions to polynomial equations, these results do not distinguish between singular and non-singular solutions. It would therefore be desirable to rephrase the conditions in Lemma~\ref{lem:Hensel} in terms of existence of $p$-adic solutions only. This is indeed possible.

We denote by $\gamma^p_d(R,m)$ the least integer $\gamma$ with the property that any system of $R$ forms of degree $d$ in $s \geq \gamma$ variables contains a $p$-adic linear space of affine dimension $m$, and define $\gamma_d^*(R,m) = \max_{p} \gamma^p_d(R,m)$. In this notation, Lemma~\ref{lem:Hensel} can be reformulated as follows.

\begin{lem}\label{lem:Qp}
    Suppose
    \begin{equation}\label{multicond}
       ms - \dim \F S_m(\B F) \geq \gamma^p_d(R,m).
    \end{equation}
    Then the system \eqref{eq:sys} possesses a non-singular $p$-adic solution, and there exists a constant $\kappa_p$ such that for any sufficiently large $\nu \in \mathbb N$ one has $\Gamma_m(p^{\nu}) \geq \kappa_p p^{\nu(ms-Rr)}$.
\end{lem}

\begin{proof}
    Let $b=\dim \F S_m(\B F)$. Recall that $\F S_m(\B F)$ really is a projective variety, so it is of projective dimension $b-1$. By a suitable version of Bertini's Theorem (e.g. \cite[Lemma 2.6]{oscar}) we can find a rational hyperplane $H_1 \subset \mathbb P_{\mathbb Q}^{ms-1}$ such that
        \[\dim_{\text{proj}} (\F S_m(\B F) \cap H_1) \leq \max\left\lbrace -1, b -2\right\rbrace \hspace{-1mm}.\]
    Thus after $k$ iterations one obtains a hyperplane section
        \[\F H_k =  H_1 \cap  \cdot \ldots \cdot \cap H_k\]
    such that the restriction of the singular locus $\F S_m(\B F) \cap \F H_k$ is of projective dimension at most $ b -k-1$. For $k=b$ this equals $-1$, so we can infer that the set $\F S_m(\B F) \cap \F H_b$ is empty as a projective variety and thus contains only the origin as an affine variety. Since all the $H_k$ are hyperplanes contained in $\mathbb P^{ms-1}_{\mathbb Q}$, the space $\F H_b$ is isomorphic to $\mathbb P^{ms-b-1}_{\mathbb Q}$ and hence of affine dimension $ms-b$. It follows that if $ms-b \geq \gamma^p_d(R,m)$ the restriction of $\F X_m(\B F)$ to $\F H_b$ is non-empty, and since $\F H_b$ does not intersect with  $\F S_m(\B F)$ we have found a non-singular point. 
\end{proof}

\section{The singular series in Birch's Theorem}

We now have the tools at hand to prove the most general version of the result. If $N_{m}(P)$ denotes the number of vectors $\B x_1, \ldots, \B x_m \in \mathbb Z^{s}$ of length at most $P$ generating a linear space contained in the complete intersection defined by $\B F = \bm 0$, we define $H_d(R,m)$ to be the least integer $H$ such that for every system $\B F$ of $R$ forms of degree $d$ the number $N_{m}(P)$ satisfies an asymptotic formula $N_{m}(P) = (c_m + o(1)) P^{ms-Rrd}$ whenever $s-\max_{\B a \ne \bm 0}\dim \F V_{\B a}(\B F)\geq H$, where $c_m$ is a non-negative constant encoding the local solution densities. In this notation our main theorem is a generalised version of Theorem~\ref{thm:asymp}. 

\begin{thm}\label{thm:asymp-m}
    Let $\B F= (F^{(1)}, \ldots, F^{(R)})$ be forms of degree $d \ge 3$ containing a non-singular real point. Furthermore, suppose that 
    \begin{align*}
      s - \max_{\B a \ne \bm 0}\dim \F V_{\B a}(\B F)\geq H_d(R,m) \quad \text{ and } \quad ms - \dim \F S_m(\B F) \geq\gamma_d^*(R,m).
    \end{align*}
    Then one has
            \begin{equation}\label{asymp-bound-m}
                N_{m}(P) =  P^{ms-Rrd} \chi_{\infty} \prod_{p \; \mathrm{ prime}} \chi_p + o\left(P^{ms-Rrd}\right) \hspace{-1mm},
            \end{equation}
    and the product of the local densities $\chi_{\infty} \prod_{p} \chi_p$ is positive.
\end{thm}
For our proof of Theorem~\ref{thm:asymp-m} we notice that by applying the arguments of \cite[\S2--3]{rd-weyl} (see also \cite[\S 7]{FRF}) to the discussion of the singularities in \cite{FRF} the condition $\max_{\B a \ne \bm 0}\dim \F V_{\B a}(\B F) \geq H_{d}(R,m)$ suffices to establish an asymptotic formula as in \eqref{asymp-bound-m} with non-negative constants $\chi_p$ and $\chi_{\infty}$. Note that for $m \geq 2$ earlier work of the author \cite[Theorem~1.1]{FRF} provides the bound
\begin{equation*}
    H_d(R,m) \leq  3 \cdot 2^{d-1}  (d-1)R(Rr+1) +1.
\end{equation*}
Furthermore, a straightforward generalisation of the arguments leading to \cite[\S2]{cubIV} shows that the real density $\chi_{\infty}$ is non-zero provided the degree $d$ is odd, so it only remains to show that if ${ms - \dim \F S_m(\B F) \geq \gamma_{d}^*(R,m)}$ the $p$-adic densities as well will be non-zero.

In \cite[eq. 6.6]{FRF} we showed that the constant $\chi_p$ can be expressed as
\begin{align}\label{p-adic}
    \chi_p &= \sum_{l=0}^{\infty}p^{-lms}\sum_{\ol x =1}^{p^l}  \dsum{\U u =1}{(\U u, p)=1}^{p^l} e \Bigg( p^{-l}\sum_{\rho=1}^R \sum_{ \B j \in J}  u^{(\rho)}_{\B j}  \Phi^{(\rho)}_{\B j}(\ol x)\Bigg),
\end{align}
 where we followed the notation from \cite{FRF} in writing $\U{u}=(u_{\B j}^{(\rho)})_{\B j \in J, 1 \leq \rho \leq R}$.
Furthermore, a standard argument (see \cite[Lemma~5.2 and Cor.]{dav} for instance) shows that the product over the $\chi_p$ converges and hence is positive as soon as every individual factor is positive. Again, it follows by a multidimensional version of the argument in \cite[Lemma~5.3]{dav} that the expression in \eqref{p-adic} can be rewritten in the shape
    \begin{equation*}
        \chi_p = \lim_{i \rightarrow \infty} p^{i(Rr-ms)}\Gamma_m(p^i),
    \end{equation*}
and by Lemma~\ref{lem:Qp} this is positive, provided the condition \eqref{multicond}
is satisfied. Replacing $\gamma^p_d(R,m)$ by $\gamma^*_d(R,m)$ now implies that the product over all $p$-adic densities is positive, which completes the proof of Theorem~\ref{thm:asymp-m}.

\end{document}